\documentclass[reqno,12pt]{amsart}
\usepackage{amssymb, verbatim}
\usepackage{graphicx}
\newtheorem{bigthm}{Theorem}
\newtheorem{thm}{Theorem}[section]
\newtheorem{lem}[thm]{Lemma}
\newtheorem{cor}[thm]{Corollary}
\newtheorem{prop}[thm]{Proposition}

\theoremstyle{definition}

\theoremstyle{remark}
\newtheorem*{rem*}{Remark}
\newtheorem{rem}[thm]{Remark}

\newtheorem*{ack*}{Acknowledgment}

\numberwithin{equation}{section}

\newcommand{\la}{\ensuremath{\lambda} }

\newcommand{\N}{\ensuremath{\mathbb N} }
\newcommand{\Om}{\ensuremath{\Omega} }

\newcommand{\C}{\ensuremath{\mathbb C} }
\newcommand{\D}{\ensuremath{\mathbb D} }
\newcommand{\whc}{\ensuremath{\widehat{\mathbb C}} }

\newcommand{\bs}{\ensuremath{\hfill _\blacksquare} }
\title{Specifying Attracting Cycles for Newton Maps of Polynomials}

\begin{document}

\author[J. Campbell]{James T. Campbell} \address[J.
Campbell]{Dept. Math. Sci.\\Dunn Hall 373\\
  University of Memphis\\Memphis, TN 38152}
\email{jcampbll@memphis.edu}

\author[J. Collins]{Jared T. Collins} \address[J. Collins]{Dept. Math. Sci.\\Dunn Hall 373\\
  University of Memphis\\Memphis, TN 38152}
\email{jtcllins@memphis.edu}

\begin{abstract}
  We show that for any set  of $n$ distinct points in the complex plane, there exists a polynomial $p$ of degree at most $n+1$ so that the corresponding Newton map, or even the relaxed Newton map, for $p$ has the given points as a super-attracting cycle. This improves the result in \cite{PlRo}, which shows how to find such a polynomial of degree $2n$. Moreover we show that in general one cannot improve upon degree $n+1$. Our methods allow us to give a simple, constructive proof of the known result that for each cycle length $n \ge 2$ and degree $d \ge 3$, there exists a polynomial of degree $d$ whose Newton map has a super-attracting cycle of length $n$. \bigskip

 \medskip 


\end{abstract}

\date{\today}

\maketitle

\section{Introduction and Statement of Results}
\label{sec:intro}

Given a complex polynomial $p$, the associated Newton map is the rational map
\[ N_p:\whc \to \whc, \quad N_p(z)=z-\frac{p(z)}{p'(z)} \, \]
and the associated relaxed Newton map is given by
\[N_{p,h}:\whc \to \whc, \quad N_{p,h}(z)=z-h\cdot\frac{p(z)}{p'(z)} \,\,  , \]
where $h$ is a complex parameter, usually taken in the disk centered at 1 of radius one. {\em Newton's method} and the {\em relaxed Newton method} are the well-known processes of iterating these maps, starting at some $z_0 \in \C$,  in search of the roots of $p$, which are attracting fixed points for the associated Newton maps.

For polynomials of degree 2, Newton's method is completely understood.  Under iteration by $N_p\, $, outside of at worst an exceptional circle in $\whc,$ every starting point in \C converges to a root of the given polynomial. For polynomials of degree three or higher, Newton's method may possess extraneous attracting (even super-attracting) cycles. Barna (\cite{b1}) seems to have been the first to notice this, and it has become a well-studied phenomenon. To illustrate how plentiful these cycles can be,  consider the one-parameter family of cubics, $p_\la(z) = z^3 + (\la -1)z - \la , \, \, \la \in \C$. Denote the relaxed Newton's map for $p_\la$ by $N_{\la, h}$.  Then there are open regions in the $\la$-plane so that for each $\la$ in these regions, there is an open subset of $h$-values for which $N_{\la, h}$ possesses extraneous attracting cycles (\cite{CuGaSu}, \cite{Kr1}). The existence of such cycles, of course, forms a barrier to using Newton's method to find the roots of the polynomial, as their basins will be open subsets of $\C$.

Many questions may be raised here. We are interested in the following. Can one specify a cycle? That is, given $n$ distinct points in the complex plane, can one construct a polynomial so that those $n$ points form an attracting cycle for the associated (relaxed) Newton's method? What is the minimal degree of such a polynomial? Conversely, fixing a degree $d$ and a cycle length $n$ (but not the cycle itself), can a polynomial of degree $d$ be constructed whose (relaxed) Newton map possesses a cycle of length $n$?

Plaza and Vergera  gave a positive answer to the first question for Newton's method. They showed that given any set $\Omega$ of $n \ge 2$ distinct points in the complex plane, there exists a polynomial $p$ of degree at most $2n$ so that $\Omega$ is a super-attracting cycle for $N_p$ (\cite{PlVe}).  Plaza and Romero (\cite{PlRo}) followed with a similar result for the relaxed Newton's method: given $\Omega$ as above, for any $h \in \D_1(1)$ there exists a polynomial $p$ of degree at most $2n$ so that $\Omega$ is a super-attracting cycle for $N_{p,h}$. (Here, $\D_1(1)$ denotes the open disk of radius 1 centered at $1$.  This disk is the largest natural region to consider for the parameter $h$ because outside of that region the roots of the polynomial become repelling points for the associated relaxed Newton map.) In  both papers they obtain their results by constructing a polynomial of degree exactly $2n$ with the prescribed super-attracting $n$-cycle.

 However their degree estimate is not minimal. Our first results show that the degree estimate may be reduced by half:

\begin{bigthm} \label{thm:a} Let $\Omega = \{z_1, z_2, \dots , z_n\}$ be any $n$ distinct points in the complex plane, $n \ge 2$. Then there exists a polynomial $p$ of degree at most $n+1$ so that $\Omega$ is a super-attracting cycle for $N_p$. \hfill \bs
\end{bigthm}

\begin{bigthm}\label{thm:b} Let $\Omega = \{z_1, z_2, \dots , z_n\}$ be any $n$ distinct points in the complex plane, $n \ge 2$. Then for any $h \in \D_1(1)$ there exists a polynomial $p$ of degree at most $n+1$ so that $\Omega$ is a super-attracting cycle for $N_{p,h}\, . \hfill \bs$
\end{bigthm}

Our methods are straightforward. For a specified cycle, we construct a homogeneous linear system so that any non-trivial solution to the system gives rise to the coefficients of the appropriate polynomial. The Rank-Nullity theorem will ensure the system has non-trivial solutions.

We also show that our estimates are sharp. For each $n \ge 2$ we demonstrate a specific collection $S_n$ of $n$ distinct points for which any polynomial whose (relaxed) Newton map has $S_n$ as a super-attracting cycle must have degree at least $n+1$.

Hurley (\cite{Hu:AO}) showed, among other things, that  for sufficiently large $k$,  there exist degree three polynomials whose Newton maps have super-attracting cycles of length $k$. Moreover he showed that such polynomials exist as perturbations by a real parameter of real polynomials, all of whose roots are real. (His results are more general, but this case may easily be gleaned from his theorems.) This gives an especially pleasing result, because Newton's method applied to real polynomials, all of whose roots are real, is known to be generally convergent (almost every (Lebesgue) initial guess lies in the attracting basin of a root). But our results allow us to prove the following:

\begin{bigthm}\label{thm:c} Let $d\geq3$ and $n \ge 2$ be natural numbers. Then there exists a polynomial $p$ of degree $d$ and a set $\Omega= \{z_1, z_2, ... , z_n\}$ of $n$ distinct points in the complex plane so that $\Omega$ is a super-attracting cycle for $N_p$. In fact the polynomial may be chosen to have real coefficients, and a cycle $\Omega$ may be found which consists entirely of real numbers.   \, \hfill \bs
\end{bigthm}

We offer this not as a new result (it is not), but because the method of proof is simple, and the polynomials found are easily described.
There is a linear system which must be examined, again to find non-trivial solutions; here however the entries in one column of the coefficient matrix in the system, after appropriate row reductions, demonstrate some recursive properties. Relatively simple techniques from one-dimensional real dynamics allow us to exploit these recursive properties, and deduce the existence of non-trivial solutions to our system.

Two interesting things come out of our proof. The first is, the process of solving  the recursive relation described in the previous paragraph, which is embedded in our (possibly quite large) linear system,  actually yields a Newton map, for a polynomial of degree d, {\em which has the same cycle we found as its own super-attracting cycle}. It is also simple to write down the corresponding polynomial, which has only one numerical parameter which must be found.

Secondly, the extraneous super-attracting $n$-cycle displays  very nice behavior. We produce real cycles of the form $\{0, 1, z_3, z_4, \dots , z_n\}$ with the property that
\[ z_1 = 0 < z_n < z_{n-1} < \cdots < z_3 < 1 = z_2 \, \, .\]

In Section \ref{sec:proof12} we prove Theorems \ref{thm:a} and \ref{thm:b}, construct the examples which show that these results are sharp, and discuss some other examples and corollaries. In Section \ref{sec:proof3} we prove Theorem \ref{thm:c}.  Because it is simpler notationally and includes all the main ideas, we prove Theorem \ref{thm:c} first in the case $d=3$ (degree of the polynomial), then follow with the proof of the general case. Demonstrations of  the properties described in the previous two paragraphs are given in Section \ref{sec:sols}. In Section \ref{sec:final} we display a pair of pictures of the basins for the Newton map found in the proof of Theorem \ref{thm:c} when $n = 5$ (the length of the extraneous cycle) and $d = 3$.

\begin{ack*}
  The authors would like to thank Paul Balister for helpful discussion regarding Theorem \ref{thm:c}. 
\end{ack*}

\section{Proofs of Theorems 1 and 2}
\label{sec:proof12}

{\bf Proof of Theorem \ref{thm:a}:} Let $\Om = \{z_1, z_2, \dots , z_n\}$ be $n$ distinct points in the complex plane, $n \ge 2$ and $f: \C \to \C$  a non-constant holomorphic function. If  \Om is a cycle for $N_f$ then it is easy to check that the following conditions are satisfied:

\begin{equation}
  \label{eq:cycle}
  f'(z_i)=\frac{f(z_i)}{z_i-z_{i+1}} \, \, , \, \, 1\leq i\leq n \, , \,  \textnormal{ where } z_{n+1}:=z_1 \, .
\end{equation}

Conversely, for any such $\Om$, if the conditions (\ref{eq:cycle}) are satisfied by a non-constant holomorphic function $f(z)$, then \Om is a cycle for $N_f$.

Suppose for the moment that $f$ is a polynomial and the distinct points in $\Omega$ form a cycle for $N_f$. Since the roots of $f$ are exactly the fixed points in \C for $N_f$,  we must have $f(z_i) \neq 0$ for every $i$, and hence $f'(z_i) \neq 0$ for every $i$. Let $\lambda_{\Omega}= \lambda$ be the multiplier for the cycle.  \Om is super-attracting for $N_f$ iff $\lambda=0$. It is straightforward to check that the condition
\begin{equation}
  \label{eq:super}
  \exists \quad z_i \in\Om\, , \, f''(z_i)=0 \quad \wedge \quad f'(z_i)\neq 0
\end{equation}
is sufficient for the cycle to be super-attracting. Again, since \Om is a cycle, the condition $f'(z_i)\neq 0$ is satisfied {\em a priori}.

Let $p(z) = a_0 + a_1z + \cdots + a_{n+1}z^{n+1}$ denote an arbitrary complex polynomial of degree at most $n+1$. We wish to find coefficients $\{a_0, a_1, \dots , a_{n+1}\}$, not all $0$,  so that (\ref{eq:cycle}) and (\ref{eq:super}) are satisfied with $f(z) = p(z)$.

Set $\vec{A} = \langle a_0,a_1, \dots , a_{n+1} \rangle$,  and for $z \in \C$ set $\vec{z} = \langle 1,z,z^2, \dots , z^{n+1}\rangle$. Then $p(z) = \vec{A} \cdot \vec{z}$, $p'(z) = \vec{A}\cdot \vec{z}\, '$, and so on for higher order derivatives.  We write $\vec{z_i}'$ for $\vec{z}\, ' |_{z = z_i}$\, and similarly for the second derivative.

 With this notation, using (\ref{eq:cycle}) and (\ref{eq:super}) we may conclude the following:
\begin{eqnarray}
  \label{eq:polycycle}
 \vec{A}\cdot((z_i-z_{i+1})\vec{z_i}'-\vec{z_i})=0  & \Leftrightarrow & N_p(z_i)=z_{i+1}   , \textnormal{ and}  \\ \label{eq:polysuper}
\vec{A}\cdot\vec{z_1}''=0 & \Rightarrow & \la = 0
\end{eqnarray}

 Now set

 \begin{equation}
   \label{eq:arrays}
   R_i = (z_i-z_{i+1})\vec{z_i}'-\vec{z_i}, \, 1\leq i \leq n, \quad R_{n+1} = \vec{z_1}'', \textnormal{ and } \, B= \left[ \begin{array}{c}         R_1\\R_2\\\vdots\\R_{n+1}\end{array} \right]\, .
 \end{equation}
The equations on the left side in each of (\ref{eq:polycycle}) and (\ref{eq:polysuper}) hold precisely when $B\vec{A}^T=\vec{0}$. Since $B$ is an $(n+1) \times (n+2)$  matrix,  the Rank-Nullity Theorem guarantees the existence of non-trivial solutions $\vec{A}$. \, \hfill \bs

\noindent \textbf{Proof of Theorem 2:} Using the above framework, set $R_i= (z_i-z_{i+1})\vec{z_i}'-h\vec{z_i}$ for $1\leq i \leq n$ and $R_{n+1} =(z_1-z_2)^2\vec{z_1}''-h(h-1)\vec{z_1}$. The rest of the argument follows as in the proof of Theorem 1.\hfill \bs

\subsection{Examples}\label{sec:examples}
    We now consider some examples. Here is the matrix  $B$ when $n = 2, 3$, and $4$:

   \[  n = 2 : \, B = \left[ \begin{array}{rrrr} -1 & -z_2 & z_1^2 - 2z_1z_2 & 2z_1^3 - 3z_1^2z_2 \\ -1 & -z_1 & z_2^2 - 2z_1z_2 & 2z_2^3 - 3z_2^2z_1 \\ 0 & 0 & 2 & 6z_1 \end{array} \right] \]

\[  n = 3 : \, B = \left[ \begin{array}{rrrrr} -1 & -z_2 & z_1^2 - 2z_1z_2 & 2z_1^3 - 3z_1^2z_2 & 3z_1^4 - 4z_1^3z_2\\ -1 & -z_3 & z_2^2 - 2z_1z_3 & 2z_2^3 - 3z_2^2z_3 & 3z_2^4 - 4z_2^3z_3\\ -1 & -z_1 & z_3^2 - 2z_1z_3 & 2z_3^3 - 3z_3^2z_1 & 3z_3^4 - 4z_3^3z_1  \\0 & 0 & 2 & 6z_1 &  12z_1^2 \end{array} \right] \]

\[  n = 4 : \, B = \left[ \begin{array}{rrrrrr} -1 & -z_2 & z_1^2 - 2z_1z_2 & 2z_1^3 - 3z_1^2z_2 & 3z_1^4 - 4z_1^3z_2 & 4z_1^5 - 5z_1^4z_2\\ -1 & -z_3 & z_2^2 - 2z_2z_3 & 2z_2^3 - 3z_2^2z_3 & 3z_2^4 - 4z_2^3z_3 & 4z_2^5 - 5z_2^4z_3\\ -1 & -z_4 & z_3^2 - 2z_3z_4 & 2z_3^3 - 3z_3^2z_3 & 3z_3^4 - 4z_3^3z_4 & 4z_3^5 - 5z_3^4z_4\\ -1 & -z_1 & z_4^2 - 2z_1z_4 & 2z_4^3 - 3z_4^2z_1 & 3z_4^4 - 4z_4^3z_1 & 4z_4^5-5z_4^4z_1 \\  0 & 0 & 2 & 6z_1 &  12z_1^2 & 20z_1^3 \end{array} \right] \]

As one can see these are beautiful matrices. In terms of solving the associated homogeneous system, however, we would prefer that they were sparser. There is one simplification which may be made, which is especially useful for small values of $n$, and will be useful later in our proof of Theorem \ref{thm:c}.  Namely, we may suppose $z_1 = 0$ and $z_2 = 1$. This is because we may always find a polynomial, of the same degree as our original, whose  Newton map has a super-attracting $n$-cycle with consecutive elements $0$ and $1$, and whose  Newton map is conjugate, via an affine transformation, to the polynomial whose Newton map had the original cycle.  Indeed, we may always find a linear fractional transformation $T$ which maps $z_1$ to 0 and $z_2$ to 1. But for any polynomial $p$, $N_p$ has a repelling fixed point at $\infty$, hence we also require $T$ to take $\infty$ to itself. These three conditions determine $T$, in fact $T(z) = (z_2-z_1)z + z_1$. A simple calculation shows that $N_{p\circ T} = T^{-1} \circ N_p \circ T$, so that $N_p$ has cycle $\{z_1, z_2, z_3\dots , z_n\}$ precisely when $N_{p \circ T}$ has cycle $\{0, 1, T(z_3), \dots , T(z_n)\}$, and the elements in this latter cycle are distinct complex numbers. The same is true in the relaxed case: here we have $N_{p\circ T, h} = T^{-1} \circ N_{p,h} \circ T$

Here are the matrices with this simplification:
\begin{eqnarray} \label{eq:simB}
  n = 2&:& B = \left[ \begin{array}{rrrr} -1 & -1 & 0 & 0 \\ -1 & 0 & 1 & 2 \\ 0 & 0 & 2 & 0 \end{array} \right] \\
 n = 3&:& B = \left[ \begin{array}{rrrrr} -1 & -1 & 0 & 0 & 0\\ -1 & -z_3 & 1 - 2z_3 & 2 - 3z_3 & 3 - 4z_3\\ -1 & 0 & z_3^2 & 2z_3^3 & 3z_3^4  \\0 & 0 & 2 & 0 &  0 \end{array} \right] \nonumber \\
  n=4&: & B = \left[ \begin{array}{rrrrrr} -1 & -1 & 0 & 0 & 0& 0\\ -1 & -z_3 & 1 - 2z_3 & 2 - 3z_3 & 3 - 4z_3 & 4 - 5z_3\\ -1 & -z_4 & z_3^2 - 2z_3z_4 & 2z_3^3 - 3z_3^2z_3 & 3z_3^4 - 4z_3^3z_4 & 4z_3^5 - 5z_3^4z_4\\ -1 & 0 & z_4^2  & 2z_4^3  & 3z_4^4  & 4z_4^5 \\  0 & 0 & 2 & 0 &  0 & 0 \end{array} \right] \nonumber
\end{eqnarray}

In particular when $n = 2$ the nullspace for $B$ is spanned by the single vector $\langle 2, -2, 0, 1 \rangle$. Thus,  any third degree polynomial with an extraneous super-attracting 2-cycle $\{0, 1\}$ produced by our method must be a multiple of $g(z) = z^3 - 2z + 2$. This is Smale's example, given in \cite{Sma:81}, to illustrate the fact that Newton's method is not generally convergent. Since for any constant $c$ and any function $f$, $f$ and $cf$ have the same Newton map, we may conclude the following:
\begin{cor}
    \label{cor:smale} If the Newton map for a degree three polynomial has  an extraneous super-attracting 2-cycle, then that Newton map is conjugate to $N_p$, where $p(z) = z^3 - 2z + 2$.
\end{cor}
\proof  Let $q$ be a degree three polynomial with an extraneous super-attracting 2-cycle $\{w_1, w_2\}$. Then $q$ must have 3 distinct roots, since it is easy to see directly that Newton maps for polynomials which have 2 or fewer roots are generally convergent, i.e., have no extraneous attracting cycles. If $T(z) = (w_2-w_1)z + w_1$ then $N_{q \circ T}$, which is conjugate to $N_q$, has the extraneous super-attracting two cycle $\{0, 1\}$. But this means that the coefficients of $p = q \circ T$ satisfy the conditions in (\ref{eq:polycycle}) for $z_1=0, z_2=1$. Since the original cycle $\{w_1, w_2\}$ is super-attracting we know that one of $q''(w_1), q''(w_2)$ are 0; renaming if necessary we may suppose that $q''(w_1) = 0$; hence the coefficients of $p = q \circ T$ also satisfy (\ref{eq:polysuper}) with $z_1 = 0$. But these conditions imply that the coefficients of $p$ must be in the nullspace of the matrix $B$ appearing in (\ref{eq:simB}), so that $p$ is a multiple of $z^3 - 2z + 2$ and hence has the same Newton map as $z^3 - 2z + 2$. But the Newton map $N_p = N_{q \circ T}$ is conjugate to the original Newton map $N_q$. \hfill \bs

\subsection{Degree Estimate is Sharp}\label{sec:sharp}

Here we construct, for each natural number $n$, a set $\Omega_n$ of $n$ distinct points so that the minimal degree polynomial whose Newton map has $\Omega_n$ as a super-attracting cycle has degree $n+1$, showing that our results are sharp.

Let $\zeta$ be a primitive $n^{th}$ root of unity (i.e. $\zeta^n =1$ and $\zeta^k \neq 1$ for $0<k<n$) and let $\Omega_n = \{\zeta, \zeta^2, \zeta^3, ..., \zeta^n = 1\}$.

Recall the matrix $B$ from the proof of Theorem \ref{thm:a}, whose entries were constructed using the elements in a specified cycle. Our first step will be to show that when $\Omega = \Omega_n$, $B$ has rank  $n+1$. Fix $\Omega_n$ and the associated matrix $B$.

Denote the $n\times n$ upper left sub-matrix of $B$ by $B_n$. That is, $B_n$ consists of the first $n$ rows of $B$, truncated after $n$ columns. The entries for $B_n$ are
\[  b_{ij}=(j-2)\zeta^{i(j-1)}-(j-1)\zeta^{i(j-2)+i+1}=\zeta^{i(j-1)}((j-2)-(j-1)\zeta), \, 1\leq i,j \leq n \, . \]
 The second factor on the right-hand side, $((j-2)-(j-1)\zeta)$, does not depend on $i$. This allows us to decompose $B_n$ in the following nice way. Let $V$ be the matrix with
 \[ v_{ij}= \zeta^{i(j-1)}, \, 1\leq i,j \leq n \, , \]
 and let $D$ be the diagonal matrix with $d_{jj}= (j-2)-(j-1)\zeta$, $1\leq j \leq n$. Then $VD=B_n$.

Since each $d_{jj} \neq 0$, $Det(D) \neq 0$. $V$ is the Vandermonde matrix generated by the cycle $\Omega = \{\zeta, \zeta^2, \zeta^3, ..., \zeta^n\}$. Thus $Det(V) = \prod_{1\leq j<i\leq n}(\zeta^i-\zeta^j) \neq 0$. Thus $Det(B_n) = Det(V)\cdot Det(D) \neq 0$, so that the rank of $B_n$ is $n$.

To determine the rank of $B$ we look at $B^T$, whose first row is $[-1,-1,...,-1,0]$ and whose $n+1^{st}$ row is \[ [(n-1)-(n)\zeta,(n-1)-(n)\zeta,...,(n-1)-(n)\zeta,(n-1)(n)\zeta^{n-2}]\, .\]
 Adding $(n-1)-(n)\zeta$ times the first row to the $n+1^{st}$ row and then dividing by an appropriate constant puts the $n+1^{st}$ row into the form $[0,0,...,0,1]$, with no change to the first $n$ rows. But the upper left $n\times n$ submatrix is $(B_n)^T$ and this matrix has rank $n$. Therefore the rank of $B^T$ is $n+1$ and in fact the reduced-row echelon form of $B^T$ will be
  \[ \left[\begin{array}{rrrrrl} 1 & 0 & 0 & \dots & 0  & \vline \, \, \, 0 \\ 0 & 1& 0 & \dots & 0 & \vline \, \, \,0 \\ 0 & 0 & 1 & \dots & 0  & \vline \, \, \, 0 \\ \vdots & \vdots & \vdots & \vdots & \vdots  & \vline \, \, \, \vdots \\ 0 & 0 & 0  & \dots & 1 & \vline \, \, \, 0 \\ \cline{1-5} 0 & 0 & 0  & \dots & 0 & \, \, \, 1    \\ 0 & 0 & 0 & \dots & 0 & \, \, \,  0 \end{array} \right] \, .\]

  Here, we have lined the boundary of the reduced-row echelon form of $B_n^T$.

  This tells us that the first $n+1$ columns of $B$ are linearly independent and since $B$ has rank $n+1$ we know the reduced-row echelon form of $B$ will be

\[ \left[\begin{array}{rrrrrrr} 1 & 0 & 0 & \dots & 0  & 0 & c_0\\ 0 & 1& 0 & \dots & 0 & 0 & c_1\\ 0 & 0 & 1 & \dots & 0  & 0 & c_2\\ \vdots & \vdots & \vdots & \vdots & \vdots  & \vdots & \vdots\\ 0 & 0 & 0  & \dots & 1  & 0 & c_{n-1}\\  0 & 0 & 0  & \dots & 0 & 1 & c_n   \end{array} \right] \, ,\]
where the $c_i$'s have values which are determined by the row-reduction process. Thus $a_{n+1}$ may be chosen as a free variable.  Choosing $a_{n+1}=0$ gives the zero polynomial, which clearly has no cycle. (Note that the zero function always satisfies (\ref{eq:cycle}), while not having \Om as a cycle for its Newton map, which is just the identity; this is the only function with this property). Any non-zero choice for $a_{n+1}$ forces a degree $n+1$ polynomial. In particular choosing $a_{n+1} = 1$ gives us the monic polynomial $p(z) = z^{n+1} - \sum_{k = 0}^n c_kz^k$ as our solution.

\section{Proof of Theorem \ref{thm:c}}
\label{sec:proof3}

\begin{rem}\label{rem:n+1}
  It follows from Theorem \ref{thm:a} that the conclusion of Theorem \ref{thm:c} is true for each pair $(d, n) = (n+1, n)$. In fact one might consider a double-induction style of proof here, but our method is general and uses neither the $(n+1, n)$ result nor induction on $(d,n)$.
\end{rem}

Recall the notation from Section \ref{sec:proof12}. If $p(z)=a_0+a_1 z+...+a_{d-1}z^{d-1}+a_d z^d$ is a complex polynomial we set $A = \langle a_0,a_1,..., a_d \rangle$ and $\vec{z} = <1,z,z^2,..., z^d>$, so that $p(z) = A \cdot \vec{z}$, $p'(z) = A\cdot\vec{z}~'$, and $p''(z) = A\cdot\vec{z}~''$. If $z_i$ is a fixed complex number we write $\vec{z_i}'$ for $\vec{z}\, ' |_{z = z_i}$\, and similarly for the second derivative.

 Given a polynomial $p$ and a set $\Omega = \{z_i\}^{n}_{i=1}$ of distinct complex numbers, the condition that $\Omega$ forms an $n$-cycle for $N_p$ is stated in (\ref{eq:polycycle}), and a sufficient condition for the cycle to be super-attracting is stated in (\ref{eq:polysuper}). Thus if we start by considering the matrix $B$ defined in (\ref{eq:arrays}), whenever  there is an $A$ with $a_d \neq 0$ satisfying $BA^T=\vec{0}$, we may conclude that $\Omega$ forms a super-attracting cycle for $N_p$, where $p$ is a polynomial of degree $d$.  For convenience here is a representation of $B$:

\[ \left[ \begin{array}{cccccc}
-1 & -z_2 & z_1^2-2z_1 z_2 & 2z_1^3-3z_1^2 z_2 & \cdots & (d-1)z_1^d-dz_1^{d-1} z_2\\
-1 & -z_3 & z_2^2-2z_2 z_3 & 2z_2^3-3z_2^2 z_3 & \cdots & (d-1)z_2^d-dz_2^{d-1} z_3\\
-1 & -z_4 & z_3^2-2z_3 z_4 & 2z_3^3-3z_3^2 z_4 & \cdots & (d-1)z_3^d-dz_3^{d-1} z_4\\
\vdots &\vdots &\vdots &\vdots &\ddots &\vdots\\
-1 & -z_{i+1} & z_i^2-2z_i z_{i+1} & 2z_i^3-3z_i^2 z_{i+1} & \cdots & (d-1)z_i^d-dz_i^{d-1} z_{i+1}\\
\vdots &\vdots &\vdots &\vdots &\ddots &\vdots\\
-1 & -z_n & z_{n-1}^2-2z_{n-1} z_n & 2z_{n-1}^3-3z_{n-1}^2 z_n & \cdots & (d-1)z_{n-1}^d-dz_{n-1}^{d-1} z_n\\
-1 & -z_1 & z_n^2-2z_n z_1 & 2z_n^3-3z_n^2 z_1& \cdots & (d-1)z_n^d-dz_n^{d-1} z_1\\
0 & 0 & 2 & 6z_1 &\cdots & (d)(d-1)z_1^{d-2}\end{array}\right] \]

We intend to argue as follows. We will view the values $z_i$ appearing in $B$ as unknown quantities, and by judicious row-reduction show that for any $n \ge 2$ and $d \ge 3$, we can always find distinct values $z_1, z_2, \dots , z_n$ so that the resulting reduced matrix, and hence $B$, has a non-trivial kernel. Moreover, considering multiplication by $B$ as a linear map from $\C^{d+1} \to \C^{n+1}$, acting on vectors $A^T$, where $A = \langle a_0,a_1,..., a_d \rangle$, we will show that the variable corresponding to $a_d$ is free, and hence may be chosen to be non-zero, in particular we may choose $a_d = 1$.  An artifact of the proof is that the $z_i$'s and $a_j$'s may be found which are real.

As noted in the proof of Theorem \ref{thm:a} we may always suppose that  $z_1=0$ and $z_2=1$. The following row operations, where the old row is on the left and its replacement is described on the right, will greatly simplify $B$:

\begin{itemize}
\item $R_1 = -R_1$
\item $R_{n+1} = \frac{R_{n+1}}{2}$
\item $R_n = R_n + R_1 - z_n^2 R_{n+1}$
\item $R_1 = R_1 - R_n$
\end{itemize}

Finally rows $R_1, R_n, R_{n+1}$ may be used to eliminate the leading terms of each of the other other rows with the row operations
\[ R_i = R_i + R_1 + z_{i+1}R_n - (z_i^2-2z_i z_{i+1})R_{n+1} \, , \, 2 \le i \le n-1 \, , \]
which  yields the following reduced form of $B$:

\[ { \Small \left[ \begin{array}{cccccc}
1 & 0 & 0 & -2z_n^3 & \cdots & -(d-1)z_n^d\\
0 & 0 & 0 & 2-3z_3-2z_n^3(1-z_3)&
\cdots & (d-1)-dz_3^{d-1} z_3-(d-1)z_n^d(1-z_3)\\
0 & 0 & 0 & 2z_3^3-3z_3^2 z_4-2z_n^3(1-z_4)&
\cdots & (d-1)z_3^d-dz_3^{d-1} z_4-(d-1)z_n^d(1-z_4)\\
\vdots &\vdots &\vdots &\vdots &\ddots &\vdots\\
0 & 0 & 0 & 2z_i^3-3z_i^2 z_{i+1}-2z_n^3(1-z_{i+1}) &
\cdots & (d-1)z_i^d-dz_i^{d-1} z_{i+1}-(d-1)z_n^d(1-z_{i+1})\\
\vdots &\vdots &\vdots &\vdots &\ddots &\vdots\\
0 & 0 & 0 & 2z_{n-1}^3-3z_{n-1}^2 z_n-2z_n^3(1-z_n)&
\cdots & (d-1)z_{n-1}^d-dz_{n-1}^{d-1} z_n-(d-1)z_n^d(1-z_n)\\
0 & 1 & 0 & 2z_n^3 & \cdots & (d-1)z_n^d\\
0 & 0 & 1 & 0 & \cdots & 0\end{array} \right] } \]

Because it is illuminating and contains all of the essential ideas of the general case,  we will present the case $d = 3$, $n \ge 3$ in the next subsection.

\subsection{The case $d=3$}
\label{sec:d=3}

\begin{rem}
  Noting Remark \ref{rem:n+1} we assume here that $n \ge 3$.
\end{rem}
When $d=3$, with the assumption that $z_1 = 0, z_2 = 1$ and after the appropriate row operations are applied, the reduced form of $B$ is an $(n+1) \times 4$ matrix of the form

$$\left[ \begin{array}{cccc}
1 & 0 & 0 & -2z_n^3\\
0 & 0 & 0 & 2-3z_3-2z_n^3+2z_n^3 z_3\\
0 & 0 & 0 & 2z_3^3-3z_3^2 z_4-2z_n^3+2z_n^3 z_4\\
\vdots &\vdots &\vdots &\vdots\\
0 & 0 & 0 & 2z_i^3-3z_i^2 z_{i+1}-2z_n^3+2z_n^3 z_{i+1}\\
\vdots &\vdots &\vdots &\vdots\\
0 & 0 & 0 & 2z_{n-2}^3-3z_{n-2}^2 z_{n-1}-2z_n^3 + 2z_n^3z_{n-1}\\
0 & 0 & 0 & 2z_{n-1}^3-3z_{n-1}^2 z_n-2z_n^3 + 2z_n^4\\
0 & 1 & 0 & 2z_n^3\\
0 & 0 & 1 & 0 \end{array}\right]$$

We will produce distinct values $\{z_1=0, z_2=1, z_3, \dots , z_n\}$ so that each of the entries in the fourth column in rows $2$ through $n-1$ are $0$ (no more row reduction). Since $z_n \neq 0$, this will imply that the family of cubic polynomials
\begin{equation}
  \label{eq:polyfam}
  p_a(z) = az^3 - 2az_n^3z + 2az_n^3, \, \, a \in \C \setminus \{0\},
\end{equation}
 has $\{0, 1, z_3, \dots , z_n\}$ as a super-attracting cycle for $N_{p_a}$. Since it will turn out that $z_n$ may be chosen to be real, and we may also choose $a$ to be real, the conclusion of Theorem \ref{thm:c} will hold when $d = 3$, for all $n$.

Simply setting the entries in the third column of row $i$ equal to $0$ ($2 \le i \le n-1$) and solving for $z_{i+1}$, one obtains the iterated relation
\begin{equation}
  \label{eq:iter} z_{i+1}=\frac{2z_n^3-2z_i^3}{2z_n^3-3z_i^2}\quad , \quad 2 \le i \le n-1 \, .
\end{equation}
Note that with $z_1 = 0, z_2=1$ and the cycle condition $z_{n+1} = z_1$, the equation appearing in (\ref{eq:iter}) is automatically satisfied for {\em any} choice of $z_n$ different from $0$ or $1$ for $i = 1$ and $i = n$.

Our proof will be complete when we find a value $z_n\neq 0, 1$ so that (\ref{eq:iter}) is satisfied for the listed values of $i$, where $z_i \neq z_j$ if $i \neq j$.  In other words, we need to find $z_n \neq 0, 1$ so that $f(z) = \frac{2z_n^3-2z^3}{2z_n^3-3z^2}$ has a cycle of length $n$ of the form $\{0, 1, z_3, \dots , z_n\}$ with $z_i \neq z_j$ for $i \neq j$. It appears that we have replaced one problem of finding a cycle with another, and indeed we have. However in the first problem we are not only finding a cycle but also a polynomial. In this second problem we have a function in hand, and this allows us to solve the problem using methods from 1-dimensional real dynamics.

Here we use notation and notions from \cite{NeAn}, where a similar situation for real polynomials is explored. Let $c \neq 0$ be a real parameter and define the family of real rational functions
\begin{equation}\label{eq:fc}
   f_c(x)=\frac{2c^3-2x^3}{2c^3-3x^2} \, , \, c \neq 0
\end{equation}
and define the sequence
\begin{equation}\label{eq:Qn}
    Q_n(c)=f_c^n(0) \, \, , \,  c \neq 0, \, \, n = 1, 2, \dots \, .
\end{equation}
In particular
\begin{eqnarray}\label{eq:Qiprops}
  Q_1(c) & \equiv & 1 \quad (c \neq 0), \textnormal{ and } \nonumber \\
  Q_{i+1}(c) &=& \frac{2c^3-2Q_i(c)^3}{2c^3-3Q_i(c)^2} \, \, .
\end{eqnarray}
Note that $Q_n(c)=0 \Leftrightarrow f_c^n(0)=0 \Leftrightarrow \, 0$ is in an $n$-cycle for $f_c$. Thus if we can find parameter values for $c$ which $Q_n(c) = 0$, we may set $z_n = c$ and we'll automatically have an $n$-cycle  for $f_c$ which includes 0. (Note that 0 goes to 1, by (\ref{eq:Qiprops}).) Then we only need to check that the minimal period is not less than $n$, and our desired result will follow.

We wish to extend the definition in (\ref{eq:Qn}) to $c = 0$. For this we do not use directly the definition involving $f_c$, rather the observation from (\ref{eq:Qiprops}) that $Q_1$ has a removable discontinuity at $c = 0$;  we extend $Q_1$ continuously to all of $\mathbb R$ by setting $Q_1(0) = 1$. The sequence $Q_n$ then defines a sequence of rational maps on $\mathbb R$, each of which is defined and continuous everywhere its denominator does not vanish.

With these conventions one easily checks that

\begin{equation}
  \label{eq:qn0}
  Q_n(0)= \left(\frac{2}{3}\right)^{n-1} , \, \, n = 1, 2, \dots \, .
\end{equation}

We wish to keep track of certain asymptotes for $Q_n$. The following observation assists us in this task:
\begin{equation}
  \label{eq:non0}
  \textnormal{For } y \in (0,\frac{3}{2}), \quad   2y^3-2Q_n(y)^3=0  \implies 2y^3-3Q_n(y)^2 \neq 0 \, , n \ge 1.
\end{equation}

Note that since $y$  is real,  $2y^3-2Q_n(y)^3=0 \Leftrightarrow  y=Q_n(y)\in (0, \frac32 )$. But if $y = Q_n(y)$ we have $2y^3-3Q_n(y)^2=2y^2(y-(3/2))$, which is non-zero in $(0, \frac32)$.

\begin{lem}
  \label{lem:alter} Let $c_n$ be the minimal positive fixed point for $Q_n$, (i.e. the minimal positive solution to $2y^3-2Q_n(y)^3=0$),  and let $d_n$ be the minimal positive solution to $2y^3-3Q_n(y)^2=0$.  Then for all natural numbers $k$,
\[ 0<...<c_k<d_k<c_{k-1}<d_{k-1}<...< 1 = c_1<d_1 = \left(\frac32\right)^{\frac13} \]
\end{lem}
\proof We employ a two-step induction. For $n=1$, $Q_1(c) = 1$ whence $y = Q_n(y)$ has minimal positive solution $c_1=1$,  while $2y^3-3(1)^2=0$ has minimal positive solution $d_1= (\frac{3}{2})^{\frac{1}{3}}>1$. Thus,  $0<1 = c_1<d_1 = \left(\frac32\right)^{\frac13}$.

Now fix $k > 1$ and suppose that $0<c_{k-1}<d_{k-1}<...<c_1<d_1$. We define
\begin{equation}
  \label{eq:gdef}
  g(y)=2y^3-3Q_k(y)^2 = 2y^3-3 (\frac{2y^3-2Q_{k-1}(y)^3} {2y^3-3Q_{k-1}(y)^2})^2 \, .
\end{equation}
By (\ref{eq:qn0}),  $g(0) =-3((\frac{2}{3})^{k-1})^2 <0$,  and by (\ref{eq:non0}),  $2(c_{k-1})^3-3Q_{k-1}(c_{k-1})^2 \neq 0$. Now
\[Q_k(c_{k-1}) = \frac{2(c_{k-1})^3-2Q_{k-1}(c_{k-1})^3}{2(c_{k-1})^3-3Q_{k-1}(c_{k-1})^2}= 0. \]
 Thus $g(c_{k-1})=2(c_{k-1})^3>0$. Moreover $g$ is rational and thus continuous where defined, which from (\ref{eq:gdef})  is everywhere except $y$-values for which $2y^3-3Q_{k-1}(y)^2=0$.  But the smallest positive solution to $2y^3-3Q_{k-1}(y)^2=0$ is $d_{k-1}$ which is outside of $(0,c_{k-1})$. Thus $g$ is defined on all of $(0, c_{k-1})$ and by the Intermediate Value Theorem there exist  solution(s) to $g(y)=0$ in the interval $(0,c_{k-1})$. The minimal such solution is $d_k$.

Next suppose that $0<d_k<c_{k-1}<d_{k-1}<...<c_1<d_1$. Let $h(y)=2y^3-2Q_k(y)^3$. By (\ref{eq:qn0}) $h(0) =-2((\frac{2}{3})^{k-1})^3 <0$. Note that $2d_k^3-3Q_k(d_k)^2 =0 \Rightarrow 2d_k^3 = 3Q_k(d_k)^2$ but $d_k<1 \Rightarrow 3Q_k(d_k)^2 <2 \Rightarrow Q_k(d_k)<\sqrt{\frac{2}{3}}$. Consider $h(d_k)= 2d_k^3-2Q_k(d_k)^3 = 3Q_k(d_k)^2 - 2Q_k(d_k)^3 = Q_k(d_k)^2[3-2Q_k(d_k)]$. Clearly $Q_k(d_k)^2>0$ and also $3-2Q_k(d_k)> 3-2\sqrt{\frac{2}{3}}>0$ thus $h(d_k)>0$. Moreover $h$ is rational and thus is continuous where defined, which is everywhere that  $2y^3-3Q_{k-1}(y)^2 \neq 0$. But the smallest positive solution to $2y^3-3Q_{k-1}(y)^2=0$ is $d_{k-1}$ which is outside of $(0,d_k)$. Hence we apply the Intermediate Value Theorem and conclude there exist solution(s) to $h(y)=0$ in the interval $(0,d_k)$. The minimal such solution is $c_k$. \bs

\begin{cor}
  \label{cor:minimal} For all  $n \ge 3$, there exists $0 < c < 1$ such that $Q_n(c)=0$,  but $Q_m(c)\neq 0$ for $m<n$.
\end{cor}
\proof By  (\ref{eq:qn0}) and minimality,  $Q_n(c_{n-1})=0$ but $Q_m(c_{n-1})\neq 0$ for $m<n$. \bs

Returning now to the proof of Theorem \ref{thm:c} in the case $d=3$,  recall that for each $n \in \N, n \ge 3$ we need to find distinct numbers $\Omega= \{0, 1, z_3, ...,z_n\}$ that satisfy the equations $$z_{i+1}=\frac{2z_n^3-2z_i^3}{2z_n^3-3z_i^2}$$
for $1\leq i\leq n$ ($z_{n+1}=z_1=0$).\\

Given $n$, choose $z_n = c_{n-1}$. Then $Q_n(c_{n-1})=0 \Leftrightarrow f_{z_n}^n(0)=0$ but recall that $f_c(x)=\frac{2c^3-2x^3}{2c^3-3x^2}$. Let $z_i = f_{z_n}^{i-1}(0)$ then $z_1 = 0$, $z_2 = \frac{2z_n^3}{2z_n^3} = 1$, $z_{n+1}=0$ and in general $$z_{i+1}=\frac{2z_n^3-2z_i^3}{2z_n^3-3z_i^2}.$$ This choice satisfies each of the equations.\\

It remains to show that the $z_i$ are distinct. Suppose that with our choice of $z_n$ we have $z_i=z_j$ for $1\leq i<j\leq n$. Then using the recursive definition above we know that: $$z_{i+1}=\frac{2z_n^3-2z_i^3}{2z_n^3-3z_i^2} \hspace{.5in} z_{j+1}=\frac{2z_n^3-2z_j^3}{2z_n^3-3z_j^2}$$
Therefore $z_i=z_j \Rightarrow z_{i+1}=z_{j+1}$ and so on.
If we choose the minimal value for $i$ so that there is a $j$ with $i<j$ and  $z_i = z_j$ then $\{0, 1, \dots, z_n\}$ will be a repeating cycle of length $j-i+1$. This tells us that our choice of $z_n$ gave a repeating cycle of length dividing $n$. In particular  we have $z_k =0$ for some $k$, $1<k<n+1$. But this forces $0 = z_k = f_{z_n}^{k-1}(0) = Q_{k-1}(c_{n-1})$,  contradicting Corollary \ref{cor:minimal}. Therefore this choice of $z_n$ yields a cycle of $n$ distinct numbers.

It is clear that since $z_n = c_{n-1}$ is real, we may choose $a$ to be real so that our resulting polynomial in (\ref{eq:polyfam}) has real coefficients, and the corresponding $n$-cycle for $N_{p_a}$ consists of real numbers.

\subsection{The proof for general $d$}
\label{sec:gen_d}

There are no real new ideas not already present in the proof for $d=3$, mostly just more bookkeeping, but we include it for the sake of completeness.

Recall that with the assumptions $z_1 = 0$ and $z_2 = 1$ the reduced form of the matrix $B$ is the following $(n+1) \times (d+1)$ matrix (for general values of $d \ge 3$ and $n \ge 3$):

\[ { \Small \left[ \begin{array}{cccccc}
1 & 0 & 0 & -2z_n^3 & \cdots & -(d-1)z_n^d\\
0 & 0 & 0 & 2-3z_3-2z_n^3(1-z_3)&
\cdots & (d-1)-dz_3^{d-1} z_3-(d-1)z_n^d(1-z_3)\\
0 & 0 & 0 & 2z_3^3-3z_3^2 z_4-2z_n^3(1-z_4)&
\cdots & (d-1)z_3^d-dz_3^{d-1} z_4-(d-1)z_n^d(1-z_4)\\
\vdots &\vdots &\vdots &\vdots &\ddots &\vdots\\
0 & 0 & 0 & 2z_i^3-3z_i^2 z_{i+1}-2z_n^3(1-z_{i+1}) &
\cdots & (d-1)z_i^d-dz_i^{d-1} z_{i+1}-(d-1)z_n^d(1-z_{i+1})\\
\vdots &\vdots &\vdots &\vdots &\ddots &\vdots\\
0 & 0 & 0 & 2z_{n-1}^3-3z_{n-1}^2 z_n-2z_n^3(1-z_n)&
\cdots & (d-1)z_{n-1}^d-dz_{n-1}^{d-1} z_n-(d-1)z_n^d(1-z_n)\\
0 & 1 & 0 & 2z_n^3 & \cdots & (d-1)z_n^d\\
0 & 0 & 1 & 0 & \cdots & 0\end{array} \right] } \]

As when $d = 3$, we proceed by demonstrating the existence of distinct values $z_1=0, z_2=1, z_3, \dots , z_n$ so that this matrix (and hence the original matrix $B$) has a non-trivial kernel. Moreover, using the notation of that previous case, we will show that the variable corresponding to $a_d$ is free, and hence may be chosen to be non-zero.

We will produce distinct values $\{z_1=0, z_2=1, z_3, \dots , z_n \neq 0\}$ so that each of the entries in the $(d+1)^{\textnormal{st}}$ column in rows $2$ through $n-1$ are $0$ (no more row reduction). Since $z_n \neq 0$, we see two things: one, that since the last column does not (and will not, even upon more row reduction) contain a ``leading 1'', the rank is at most $d$; and second, non-trivial solutions to the homogeneous system may be found by choosing  $a_d$ (the coefficient on $z^d$ in our polynomial) to be non-zero. Taking any such solution (there may be other free variables) will give us a polynomial of degree $d$ with  $\{0, 1, z_3, \dots , z_n\}$ as a super-attracting cycle for $N_{p_a}$. Since it will turn out that $z_n$ may be chosen to be real, we may also choose $a$ to be real and the conclusion of Theorem \ref{thm:c} will hold.

Simply setting the entries in the $(d+1)^{\textnormal{st}}$ column of row $i$ equal to $0$ ($2 \le i \le n-1$) and solving for $z_{i+1}$, one obtains the iterated relation
\begin{equation}
  \label{eq:iterd} z_{i+1}=\frac{(d-1)z_n^d-(d-1)z_i^d}{(d-1)z_n^d-dz_i^{d-1}} \quad , \quad 2 \le i \le n-1 \, .
\end{equation}
Note that with $z_1 = 0, z_2=1$ and the cycle condition $z_{n+1} = z_1$, the equation in (\ref{eq:iterd}) is automatically satisfied for {\em any} choice of $z_n$ different from $0$ or $1$, for $i = 1$ and $i = n$.

Let $c \neq 0$ be a real parameter and define the family of real rational functions
\begin{equation}\label{eq:fcd}
   f_c(x)=\frac{(d-1)c^d-(d-1)x^d}{(d-1)c^d-dx^{d-1}} = \frac{c^d - x^d}{c^d- (\frac{d}{d-1})x^{d-1}}\,
\end{equation}
and define the sequence
\begin{equation}\label{eq:Qnd}
   Q_n(c)=f_c^n(0) \, \, , \, c \neq 0 \, , \,  n = 1, 2, \dots \, .
\end{equation}
In particular
\begin{eqnarray}\label{eq:Qid}
  Q_1(c)&\equiv & 1 \, (c \neq 0), \textnormal{ and } \nonumber \\
  Q_{i+1}(c) &=& \frac{c^d-Q_i(c)^d}{c^d-(\frac{d}{d-1})Q_i(c)^{d-1}} \, \, .
\end{eqnarray}
Note that $Q_n(c)=0 \Leftrightarrow f_c^n(0)=0 \Leftrightarrow f_c$ has an n-cycle including zero.

We wish to extend the definition in (\ref{eq:Qnd}) to $c = 0$. For this we do not use directly the definition involving $f_c$, rather the observation from (\ref{eq:Qid}) that $Q_1$ has a removable discontinuity at $c=0$; we extend $Q_1$ continuously to all of $\mathbb R$ by setting $Q_1(0) = 1$. The sequence $Q_n$ then defines a sequence of rational maps on $\mathbb R$, each of which is defined and continuous everywhere its denominator does not vanish.

With these conventions one easily checks that

\begin{equation}
  \label{eq:qn0d}
  Q_n(0)= \left(\frac{d-1}{d}\right)^{n-1} , \, \, n = 1, 2, \dots \, .
\end{equation}
We wish to keep track of certain asymptotes for $Q_n$. The following observation helps us with this task:

\begin{eqnarray} \label{eq:non0d}
  \textnormal{For } y  \in  (0,\frac{d}{d-1}) & \textnormal{ and all} & n \ge 1,    \nonumber \\
  y^d-Q_n(y)^d=0 & \implies &  y^d- (\frac{d}{d-1}) Q_i(y)^{d-1} \neq 0 \,.
\end{eqnarray}

 Observe that since $y$ is positive,  $y^d-Q_n(y)^d=0 \Rightarrow y=Q_n(y) \in (0,\frac{d}{d-1})$. But if $y = Q_n(y)$ we have $(d-1)y^d-dy^{d-1} =y^{d-1}((d-1)y-d)$, which is non-zero in $(0,\frac{d}{d-1})$.

\begin{lem}
  \label{lem:alterd} For a fixed $d$ let $c_n$ be the minimal positive solution to \[y^d-Q_n(y)^d=0 \, , \]
  and let $b_n$ be the minimal positive solution to
  \[(d-1)y^d-dQ_n(y)^{d-1} = 0 \, . \]
  Then for all natural numbers $k$,
  \begin{equation}\label{eq:alterd}
   0<...<c_k<b_k<c_{k-1}<b_{k-1}<...<c_1=1 <b_1 = (\frac{d}{d-1})^{\frac{1}{d}}\, .
  \end{equation}
   \end{lem}

\proof
For $n=1$, $Q_1(c) =1$ thus $y^d-(1)^d=0$ has minimal positive solution $c_1=1$ while $(d-1)y^d-d(1)^{d-1}=0$ has minimal positive solution $b_1= (\frac{d}{d-1})^{\frac{1}{d}}>1$, so that the $0 < c_1 = 1 < b_1 = (\frac{d}{d-1})^{\frac{1}{d}}$.

Suppose that $0<c_{k-1}<b_{k-1}<...<c_1<b_1$. Let $g(y)=(d-1)y^d-dQ_k(y)^{d-1}$. By (\ref{eq:qn0d}), $g(0) =-d((\frac{d-1}{d})^{k-1})^{d-1} <0$ and by (\ref{eq:non0d}) $(d-1)(c_{k-1})^d-dQ_{k-1}(c_{k-1})^{d-1} \neq 0$.

 \[Q_k(c_{k-1}) = \frac{(d-1)(c_{k-1})^d-(d-1)Q_{k-1}(c_{k-1})^d}{(d-1)(c_{k-1})^d-dQ_{k-1}(c_{k-1})^{d-1}}= 0.\]
Thus $g(c_{k-1})=(d-1)(c_{k-1})^d>0$. Moreover $g$ is rational and thus is continuous where defined. $g(y) = (d-1)y^d-d (\frac{(d-1)y^d-(d-1)Q_{k-1}(y)^d} {(d-1)y^d-dQ_{k-1}(y)^{d-1}})^{d-1}$ but the smallest solution to $(d-1)y^d-dQ_{k-1}(y)^{d-1}=0$ is $b_{k-1}$ which is outside of $(0,c_{k-1})$. Thus by the Intermediate Value Theorem there exist solution(s) to $g(y)=0$ in the interval $(0,c_{k-1})$.  Call the minimal solution $b_k$.\\

Next suppose that $0<b_k<c_{k-1}<b_{k-1}<...<c_1<b_1$. Let $h(y)=(d-1)y^d- (d-1)Q_k(y)^d$. By (\ref{eq:qn0d}),  $h(0) =-(d-1)((\frac{d-1}{d})^{k-1})^d <0$. Note that $(d-1)b_k^d-dQ_k(b_k)^{d-1} =0 \Rightarrow (d-1)b_k^d = dQ_k(b_k)^{d-1}$
but $0<b_k\leq(\frac{d}{d-1})^{\frac{1}{d}} \Rightarrow 0< dQ_k(b_k)^{d-1} \leq1 \Rightarrow 0<Q_k(b_k)<1$. Consider $h(b_k)= (d-1)b_k^d-(d-1)Q_k(b_k)^d = dQ_k(b_k)^{d-1} - (d-1)Q_k(b_k)^d = Q_k(b_k)^{d-1}[d-(d-1)Q_k(b_k)] >0$.  Moreover $h$ is rational and thus is continuous where defined. $h(y) = (d-1)y^d-(d-1) (\frac{(d-1)y^d- (d-1)Q_{k-1}(y)^d} {(d-1)y^d-dQ_{k-1}(y)^{d-1}})^d$ but the smallest solution to $(d-1)y^d-dQ_{k-1}(y)^{d-1}=0$ is $b_{k-1}$ which is outside of $(0,b_k)$. Again by the Intermediate Value Theorem there exist  solution(s) to $h(y)=0$ in the interval $(0,b_k)$, call the minimal solution $c_k$. \bs

As in the case $d=3$ we now have
\begin{cor}
  \label{cor:minimald} Let $d \in \mathbb N$, $d \ge 3$. For all  $n \ge 3$, there exists $0 < c < 1$ such that $Q_n(c)=0$,  but $Q_m(c)\neq 0$ for $m<n$. \bs
\end{cor}

To finish the proof of Theorem \ref{thm:c} for general $d$, proceed as in the case $d=3$: for a fixed $d$ and $n$ choose $z_n = c_{n-1}$ found in Lemma \ref{lem:alterd}. Then $Q_n(c_{n-1})=0 \Leftrightarrow f_{z_n}^n(0)=0$. Recall that $f_c(x)=\frac{(d-1)c^d-(d-1)x^d}{(d-1)c^d-dx^{d-1}}$. Let $z_i = f_{z_n}^{i-1}(0)$;  then $z_1 = 0$, $z_2 = \frac{(d-1)z_n^d}{(d-1)z_n^d} = 1$, $z_{n+1}=0$ and in general $$z_{i+1}=\frac{(d-1)z_n^d-(d-1)z_i^d}{(d-1)z_n^d-dz_i^{d-1}}.$$ This choice satisfies each of the equations.

It remains to see that the $z_i$ are distinct. Suppose that a choice of $z_n$ forces $z_i=z_j$ for $1\leq i<j\leq n$. Then using the recursive definition above we know that: $$z_{i+1}=\frac{(d-1)z_n^d-(d-1)z_i^d}{(d-1)z_n^d-dz_i^{d-1}} \hspace{.5in} z_{j+1}=\frac{(d-1)z_n^d-(d-1)z_j^d}{(d-1)z_n^d-dz_j^{d-1}}$$
Therefore $z_i=z_j \Rightarrow z_{i+1}=z_{j+1}$ and so on. If we choose the minimal value for $i$ so that there is a $j$ with $i<j$ and  $z_i = z_j$ then $\{0, 1, \dots , z_n\}$ will be a repeating cycle of length $j-i+1$. This tells us that the only bad choices of $z_n$ are those that lead to repeating cycles of length dividing $n$. In particular those cycles with $z_j =0$ for $1<j<n+1$. But this forces $0 = z_j = f_{z_n}^{j-1}(0) = Q_{j-1}(c_{n-1})$ contradicting Corollary \ref{cor:minimald} . Therefore the $z_i$ are distinct.

Once again it is clear that since $z_n = c_{n-1}$ is real, and we may choose $a_d$ to be real, our resulting polynomial in has real coefficients, and the corresponding $n$-cycle for $N_{p_a}$ consists of real numbers. This concludes the proof of Theorem \ref{thm:c}. \bs

\section{Properties of the Solutions}
\label{sec:sols}
The rational map
\begin{equation}\label{eq:fcdef}
   f_c(z)=\frac{(d-1)c^d-(d-1)z^d}{(d-1)c^d-dz^{d-1}} = \frac{c^d-z^d}{c^d-\left(\frac{d}{d-1}\right)z^{d-1}}
\end{equation}
examined in the proof of Theorem \ref{thm:c} is the Newton map for the polynomial
\[ p_c(z) = z^d - (d-1)c^dz + (d-1)c^d.  \]
As part of the proof, for fixed $d \ge 3$,  we found that if $c = c_{n-1}$, where $c_k$ is defined for general $k$ in Lemma \ref{lem:alterd}, then $f_c(z) = N_{p_c}(z)$ has a super-attracting $n$-cycle $\{0, 1, z_3, z_4, \dots , z_n = c\}$. Thus, up to the question of actually finding $c_{n-1}$, we have a simple description of a degree $d$ polynomial with a super-attracting $n$-cycle.

It turns out that the following nice property holds for the cycle:

\begin{prop}
  \label{prop:cycleorder} Fix $d\ge 3$. Then for all $n \ge 3$ the cycle $\Omega = \{0, 1, z_3 \dots , z_n\}$ found in Theorem \ref{thm:c} satisfies
  \begin{equation}
    \label{eq:order}
    z_1 = 0 < z_n < \dots < z_{i+1} < z_i < \dots < z_3 < 1 = z_2 \, .
  \end{equation}
\end{prop}
\proof Fix $d\ge 3$. Because we will be comparing cycle elements from cycles of different lengths, we introduce double subscripts to keep track of elements from cycles of different lengths. Also, to ease notation, we re-number the sequence $\{c_k\}$ by $k \to k+1$. Hence we will write

\begin{equation}\label{eq:recycle}
    0 = z_{1, n}, \, 1  = z_{2,n}= f_{c_n}(0), \, z_{3, n} = f_{c_n}^2(0), \dots , \, z_{n, n} =  f_{c_n}^{n-1}(0) = c_n
\end{equation}
for the cycle of length $n$. With this notation, the following holds (e.g. (\ref{eq:alterd})):
\begin{equation}\label{eq:znn}
    0 < \dots < z_{k+1, k+1} < z_{k, k} < \dots < z_{3, 3} < z_{2, 2} = 1.
\end{equation}

We will be considering the maps $f_c(x)$ only for $0 < c < 1$. For $c$ in this range, the following useful properties of $f_c$ are easily established:
\begin{eqnarray}
  \label{eq:asymp} & &f_c(x) \textnormal{ has a vertical asymptote at } x_0= \left(\frac{d-1}{d}\right)c^{\frac{d}{d-1}}\, < \, c \, . \\
  \label{eq:incx} & &f_c(x) \textnormal{ is an increasing function of } x, x > 0, \, x \neq x_0   \\
  \label{eq:decc} & &\textnormal{If } 0 < c < c' \le x < \frac{d}{d-1}, \, \textnormal{ then } f_{c'}(x) < f_{c}(x)
\end{eqnarray}
To finish the proof of the Proposition, we need the following lemmas:
\begin{lem}\label{lem:moveright}
\textnormal{For all } $k \ge 3$,
\begin{equation}\label{eq:moveright}
    0 < c_k = z_{k, k} < z_{k, k+1} < z_{k, k+2} < \cdots \, 1\, .
\end{equation}
\textnormal{That is, as the cycle length increases, the $k^{\textnormal{th}}$ element of the cycle moves to the right.}
\end{lem}
 \noindent {\bf Proof of Lemma:} We induct on $k$. Note that $0 < c_k$ for all $k$. Set $k = 3$ and recall that $z_{3,j}=f_{c_j}(1)< 1$. From Lemma \ref{lem:alterd}, we know that the $c_j$'s are decreasing (as j increases), so from (\ref{eq:decc}) we may conclude that when $i < j$ we have $  z_{3,i}=f_{c_i}(1) < f_{c_j}(1)= z_{3,j}< 1$, that is, (\ref{eq:moveright}) holds for $k = 3$.

Suppose that for $k-1$, $0 < z_{k-1,k-1}<z_{k-1,k}<z_{k-1,k+1}<\dots < 1$. Suppose $j > i \ge k$.
  Since $f_{c_j}(x)$ is a increasing function of $x$ we have $z_{k-1,i}<z_{k-1,j} \Rightarrow f_{c_j}(z_{k-1,i})<f_{c_j}(z_{k-1,j}).$ Again applying (\ref{eq:decc}) we have
  \[z_{k,i} = f_{c_i}(z_{k-1,i})<f_{c_j}(z_{k-1,i})<f_{c_j}(z_{k-1,j}) = z_{k,j}< f_{c_j}(1) < 1 \]
  and since $0 < c_k = z_{k,k}$, the Lemma is proved. \bs

 \begin{cor}\label{cor:goodorder}
 \begin{equation}\label{eq:lastobs}
    \textnormal{For } 3 \le i < j \, , \, c_j = z_{j,j} < z_{i,j} \, .
 \end{equation}
 \end{cor}
That is, as you keep the cycle length fixed at $n=j$, the elements of the cycle (except for $z_{1, j} = 0$) don't move to the left of $c_j$.
\proof Indeed, when $i < j$, $c_j = z_{j,j} < z_{i,i} < z_{i,j}$. \bs

For the last step in the proof of the Proposition, we recall that $d \ge 3$ is fixed, and we fix $n \ge 4$. (The Proposition is true when $n = 3$). We dispense with the double subscripts and label our cycle as $z_1 = 0, z_2, = 1, z_3, \dots , z_n = c_n$. With this notation Corollary \ref{cor:goodorder} tells us that for all $3 \le i < n$, $z_n < z_i$. This is enough to give us the required order (\ref{eq:order}). We have the following string of elementary implications:

$$(d-1)z_n^d (1-z_i)+z_i^d>0$$
$$\Rightarrow (d-1)z_n^d + z_i^d > (d-1)z_n^d z_i$$
$$\Rightarrow (d-1)z_n^d -(d-1)z_i^d > (d-1)z_n^d z_i - dz_i^d$$
$$\Rightarrow (d-1)z_n^d -(d-1)z_i^d > [(d-1)z_n^d - dz_i^{d-1}]z_i$$
$$\Rightarrow z_{i+1} < z_i, $$

and as we know that $z_n < z_{i+1}$, the Proposition follows.  \bs

\section{Graphics} \label{sec:final}

Here are two pictures of the attracting basin for the Newton map $f_c$ when $d = 3$ and $n = 5$. In this case, $c_5 \sim 0.245772$ and $f_c$ is the Newton map for the polynomial $p_c(z) = z^d - (d-1)c^dz + (d-1)c^d = z^3 -0.0296913z + .0296913$.

The first picture is centered at the origin, in a square window of sidelength approximately 1.5. As might be expected, because the roots of $p_c$ are spaced fairly evenly and the 5-cycle is constricted to the interval $[0,1]$, the picture has a lot of symmetry:
\begin{figure}[h]
\caption[bigbasin]{Basins of attraction, d=3, n=5}
\includegraphics[scale=0.95]{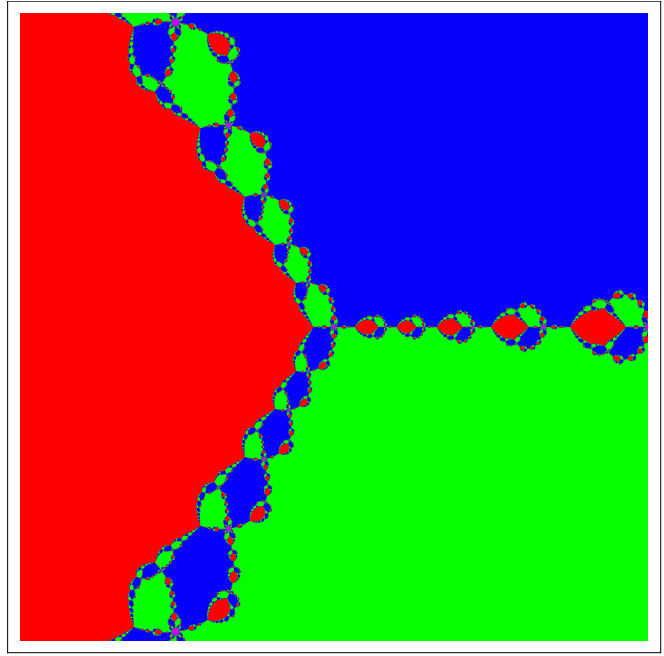}
\end{figure}

\pagebreak

Here is a closeup in the window $[-.05, .255] \times [-.125, .125]$. The basins for the attracting cycle appear on the $x$-axis at the base of the flower-like figure, colored purple, and appear to be disks.

\begin{figure}[h]
\caption[smallbasin]{Close-up of Basins, near the origin}
\includegraphics[scale=0.95]{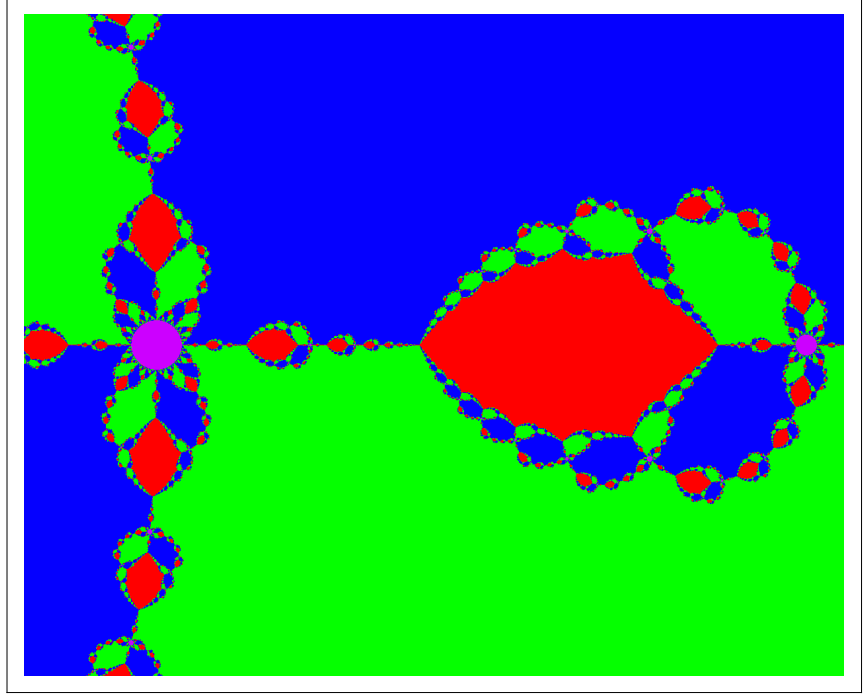}
\end{figure}

\bibliographystyle{amsalpha}
\bibliography{newton}

\end{document}